\documentclass[reqno]{article}

\usepackage{graphicx,graphics,subfigure,enumerate,color,epsfig,psfrag}

\usepackage{amsmath,amssymb,amsfonts,bm}
\usepackage{empheq}
\usepackage{comment}

\usepackage{hyperref}
\usepackage{setspace}
 \evensidemargin=\oddsidemargin
\newdimen\dummy
\dummy=\oddsidemargin 
\newcommand{\com}[1]{{\color{blue}{#1}}}

\newcommand{\Real}{\mathbf R}
\newcommand{\eps}{\varepsilon}

\newcommand{\Bb}{\mathbf{b}}

\newcommand{\vb}{\mathbf v}
\newcommand{\wb}{\mathbf w}
\newcommand{\xb}{\mathbf x}
\newcommand{\ub}{\mathbf u}
\newcommand{\Fb}{\mathbf F}
\newcommand{\Jm}{\mathbb{J}}
\newcommand{\Id}{\mathbb{I}}

\newcommand{\gad}{gentlest ascent dynamics }

\begin{document}
\begin{center}
{\bf \large The Gentlest Ascent Dynamics}

Weinan E

Department of Mathematics and PACM,
Princeton University

Xiang Zhou

Division of Applied Mathematics,
Brown University

\end{center}

\begin{abstract}
Dynamical systems that describe the escape from the basins of attraction of
stable invariant sets  are presented and analyzed.
It is shown that the stable fixed points of such dynamical systems are the
index-1 saddle points.  Generalizations to high index saddle
points are discussed.  Both gradient and non-gradient systems are considered.
Preliminary results on the nature of the dynamical behavior are presented.
\end{abstract}

\section{The gentlest ascent dynamics}

Given an energy function  $V$ on $\Real^n$, 
the simplest form of the steepest decent dynamics (SDD) associated with $V$ is
\begin{equation}
\label{SDD-grad}
 \dot{\xb} = - \nabla V(\xb).
\end{equation}
It is easy to see that if $\xb(\cdot)$ is a solution to \eqref{SDD-grad}, then
$V(\xb(t))$ is a decreasing function of $t$.  Furthermore, the stable fixed points of
the dynamics \eqref{SDD-grad} are the local minima of $V$.
Each local minimum has an associated basin of attraction which consists of all the initial
conditions from which the dynamics described by \eqref{SDD-grad} converges to that local
minimum as time goes to infinity. 
For \eqref{SDD-grad}, these are simply the potential wells of $V$.
The basins of attraction are separated by separatrices, on
which the dynamics converges to saddle points.

We are interested in the opposite dynamics: 
The dynamics of escaping a basin of attraction.
The most naive suggestion is to just reverse the sign in \eqref{SDD-grad}, 
the dynamics would then find the local maxima of $V$ instead.  
This is not what we are interested in.
We are interested in the gentlest
way in which the dynamics climb out of the basin of attraction.
Intuitively, it is clear that what we need is a dynamics that converges to the index-1
saddle points of $V$.
Such a problem is of general interest to the study of noise-induced transition between
metastable states \cite{String2002,RMP1990} :
Under the influence of small noise, with high probability,
the escape pathway  has to go through the neighborhood of a saddle point \cite{FW1998}.

The following dynamics serves the purpose:
\begin{subequations}
\begin{center}
\begin{empheq}[left=\empheqlbrace]{align}
 \dot{\xb}  &= - \nabla V(\xb) + 2\frac{(\nabla V, \vb)}{(\vb, \vb)} \vb,  \\
 \dot{\vb}  &= - \nabla^2 V(\xb)\vb + \frac{(\vb, \nabla^2 V \vb)}{(\vb, \vb)} \vb  .
 \end{empheq} 
\end{center}
\label{GAD-grad}
\end{subequations}

We will show later that the stable fixed points of this dynamics are precisely the index-1
saddle points of $V$ and the unstable directions of $V$ at the saddle points.
Intuitively the idea is quite simple:  The second equation in \eqref{GAD-grad}
attempts to find the direction that corresponds to the smallest eigenvalue of 
$ \nabla ^2 V$,
and the last term in the first equation makes this direction an ascent direction.

This consideration is not limited to the so-called ``gradient systems'' such as
\eqref{SDD-grad}. It can be extended to  non-gradient systems.
Consider the following dynamical system:
\begin{equation}
\label{SDD-non-grad}
 \dot{\xb} = \Fb (\xb).
\end{equation}
We can also speak about the stable invariant sets of this system, and escaping
basins of attraction of the stable  invariant sets.
In particular, we can also think about finding index-1 saddle points,
though in this case, there is no guarantee that under the influence of small
noise, escaping the basin of attraction has to proceed via saddle points\cite{Xiangthesis}.

For non-gradient systems, \eqref{GAD-grad} has to be modified to
\begin{subequations}
\begin{center}
\begin{empheq}[left=\empheqlbrace]{align}
  \dot{\xb} &= \Fb(\xb) -2 \frac{(\Fb(\xb), \wb )}{ (\wb,\vb)}\vb,\label{GAD-non-grad-x}\\
 \dot{\vb} &= (\nabla \Fb(\xb))\vb - \alpha(\vb)\vb, \label{GAD-non-grad-v}\\
 \dot{\wb} &= (\nabla \Fb(\xb))^T \wb - \beta(\vb,\wb)\wb. \label{GAD-non-grad-w}
\end{empheq} 
\end{center}
\label{GAD-non-grad}
\end{subequations}
Here two directional vectors $\vb$ and $\wb$ are needed in order to follow
both the right and left eigenvectors of the Jacobian.
Given the matrix $\nabla \Fb(\xb)$, two scalar valued functions $\alpha$ and $\beta$  are defined by 
\begin{subequations}
\begin{center}
\begin{empheq}[left=\empheqlbrace]{align}
\alpha(\vb) &=(\vb, (\nabla \Fb(\xb))\vb), \\
\beta(\vb,\wb)&= 2(\wb, (\nabla \Fb(\xb))\vb)-\alpha(\vb).
 \end{empheq} 
\end{center}
\label{eqn:alpha-beta}
\end{subequations}
We have taken and we will take the normalization such that $(\vb,\vb)=1$ and $(\wb, \vb) = 1$.
They are  to keep the normalization such that $(\vb,\vb)=1$ and $(\wb, \vb) = 1$. 
This normalization is preserved by the dynamics as long as it holds initially.
Thus, the first equation in \eqref{GAD-non-grad} actually is equivalent to   $\dot{\xb} = \Fb(\xb) -2 {(\Fb(\xb), \wb )} \vb$.
(Of course, one can enforce other types of normalization condition, such as the symmetric one:
$(\vb,\vb)=(\wb,\wb)$ and $(\wb, \vb) = 1$, and define new expressions of $\alpha$ and $\beta$ accordingly.)
  In the case of gradient flows, we can take $\wb=\vb$ and \eqref{GAD-non-grad}
reduces to \eqref{GAD-grad}.

\begin{figure}[htbp!]
\begin{center}
\scalebox{0.8}{ \input{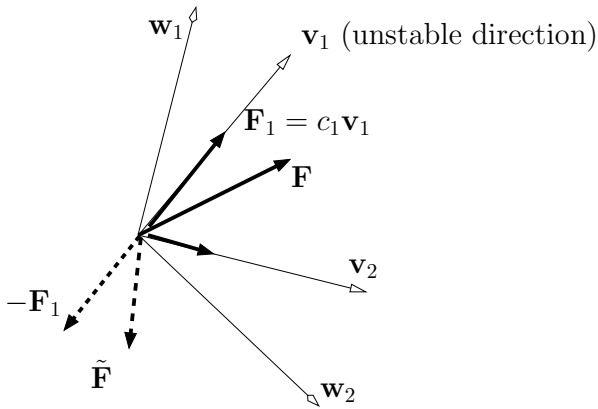} }\ \caption{Illustration of the 
gentlest ascent dynamics.
$\Fb$ is the force of the original dynamics and $\tilde{\Fb}$ is the
force of the gentlest ascent dynamics. $\vb_1$ and $\vb_2$ represent the
unstable and stable right eigenvectors, respectively; $\wb_1$ and
$\wb_2$ are the corresponding left eigenvectors. Note that $\wb_1\perp
\vb_2$ and $\wb_2\perp \vb_1$. $\Fb$ has the decomposition
$\Fb=\Fb_1+\Fb_2=c_1\vb_1+c_2\vb_2$ where the coefficient $ c_1 = (\Fb,
\wb_1)/(\vb_1, \wb_1)$. Thus,   $\tilde{\Fb}:=-\Fb_1+\Fb_2=\Fb-2\Fb_1=\Fb-2c_1\vb_1$. }
\label{fig:EF_illustrate}
\end{center}
\end{figure}

We call this the {\it gentlest ascent dynamics}, abbreviated GAD.  
It has its origin in some of the numerical
techniques proposed for finding saddle points.  For example, there is indeed a numerical
algorithm proposed by Crippen and Scheraga called  the
``gentlest ascent method'' \cite{Crippen1971}.  The main idea
is similar to that of GAD, namely to find the right direction, the direction of the
eigenvector corresponding to the smallest eigenvalue and making that an ascent
direction.
But the details of the gentlest ascent method seem to be quite a bit more complex.
The ``eigenvector following method'' proposed in literature, for example, \cite{cerjan1981,Energylanscapes},
is based on a very similar idea.
There at each step, one finds the eigenvectors of the Hessian matrix of the potential.
Also closely related is the ``dimer method'' in which two states connected by a small line
segment are evolved simultaneously in order to find the saddle point \cite{Dimer1999}.
One advantage of the dimer method is that it avoids computing the Hessian of the potential.
From the viewpoint of our GAD,
the spirit of ``dimer method'' is equivalent to use central  difference scheme to numerically calculate the matrix-vector multiplication 
in GAD \eqref{GAD-non-grad} and \eqref{eqn:alpha-beta} by writing
$(\nabla \Fb(\xb)) \Bb =  \frac{d }{d \eps} \Fb(\xb+\eps\Bb)|_{\eps=0}\approx \frac{1}{2\eps} (\Fb(\xb+\eps\Bb)-\Fb(\xb-\eps\Bb))$ for any vector $\Bb$.

We believe that as a dynamical system, the continuous formulation embodied in
\eqref{GAD-grad} and \eqref{GAD-non-grad} has its own interest.
We will demonstrate some of these interesting aspects in this note.

\noindent
{\bf Proposition.}
Assume that the vector field $\Fb$ is $C^3(\Real^{n})$. 
\begin{enumerate}[(a)]
\item  If $(\xb_{*},\vb_{*},\wb_{*})$ is a fixed point of the \gad \eqref{GAD-non-grad} and  
$\vb_{*},\wb_{*}$ are normalized such that $\vb_{*}^{T}\vb_{*}=\vb_{*}^{T}\wb_{*}=1$, then 
$\vb_{*}$ and $\wb_{*}$ are the   right and left eigenvectors , respectively, of $\nabla \Fb(\xb_{*})$ corresponding to one eigenvalue $\lambda_{*}$, i.e.,
$$(\nabla \Fb(\xb_{*}))\vb_{*}=\lambda_{*} \vb_{*} , \quad  (\nabla \Fb(\xb_{*}))^{T}\wb_{*}=\lambda_{*} \wb_{*},$$
and $\xb_{*}$ is a fixed point of the original dynamics system, i.e., $\Fb(\xb_{*})=\mathbf{0}$.

\item Let $\xb_s$ be a fixed
point of the original dynamical system  $\dot{\xb}=\Fb(\xb)$. 
If the Jacobian matrix $\Jm(\xb_s)=\nabla \Fb(\xb_{s})$ has $n$ 
distinct real eigenvalues
$\lambda_1, \lambda_2, \cdots, \lambda_n$ and  $n$ linearly independent right and left
eigenvectors, denoted by $\vb_i$ and $\wb_i$  correspondingly,  i.e.,
$$\Jm(\xb_{s})\vb_{i}=\lambda_{i}\vb_{i},\qquad \Jm(\xb_{s})^{T}\wb_{i}=\lambda_{i}\wb_{i},\quad i=1,\cdots, n$$
and in addition, we impose 
the normalization condition $\vb_{i}^{T}\vb_{i}=\wb_{i}^T  \vb_i=1,\ \forall i$ ,  
then 
for all $i=1,\cdots,n$,  $(\xb_s,\vb_i,\wb_i)$ is a fixed point of the \gad \eqref{GAD-non-grad}.
Furthermore, among these $n$ fixed points, there exists
one fixed point $(\xb_s,\vb_{i'},\wb_{i'})$ which is linearly stable
 if and only if $\xb_s$ is an index-$1$ saddle point  of the original dynamical system
 $\dot{\xb}=\Fb(\xb)$ and the eigenvalue $\lambda_{i'}$ corresponding to $\vb_{i'}$, $\wb_{i'}$
 is the only positive eigenvalue of $\Jm(\xb_s)$.

\end{enumerate}

\hspace{-6mm}
{\bf Proof.} 

(a) Under the given condition, it is obvious that $
(\nabla \Fb(\xb_{*}))\vb_{*} = \alpha(\vb_{*})\vb_{*}$ and 
$(\nabla \Fb(\xb_{*}))^T \wb_{*} = \beta(\vb_{*},\wb_{*})\wb_{*}$.
By definition and other conditions,
$
\beta(\vb_{*},\wb_{*}) = 2\wb_{*}^{T} (\nabla \Fb(\xb_{*}))\vb_{*}-\alpha(\vb_{*})=
2\wb_{*}^{T} (\alpha(\xb_{*}))\vb_{*}-\alpha(\vb_{*})=\alpha(\vb_{*})
$. Therefore,  $\vb_{*}$ and $\wb_{*}$ share the same eigenvalue $\lambda_{*}=\alpha(\vb_{*})=\beta(\vb_{*},\wb_{*})$.
From the fixed point condition 
$\Fb(\xb_{*})-2 (\wb_{*}^{T}\Fb(\xb_{*})) \vb_{*}=\mathbf{0}$, 
we take the inner product of this equation with $\wb_{*}$ to get   $\wb_{*}^{T}\Fb(\xb_{*}) = 2\wb_{*}^{T}\Fb(\xb_{*}) $.
So $\wb_{*}^{T}\Fb(\xb_{*})=0$ and in consequence, the conclusion $\Fb(\xb_{*})=\mathbf{0}$ holds from the 
fixed point condition $\Fb(\xb_{*})-2 (\wb_{*}^{T}\Fb(\xb_{*})) \vb_{*}=\mathbf{0}$ again.

(b)
It is obvious that for all $i$,
$(\xb_s,\vb_i,\wb_i)$ is a fixed point of the \gad \eqref{GAD-non-grad} 
by the definition of $\vb_i$ and $\wb_i$. It is going to be shown that
we can explicitly write down the eigenvalues and eigenvectors of GAD at any fixed point
$(\xb_s,\vb_i,\wb_i)$.

Let $\Jm(\xb)= \nabla \Fb(\xb)$.
The Jacobian matrix of the \gad \eqref{GAD-non-grad}  has the following expression: \
$ \tilde{\Jm}(\xb,\vb,\wb)  =$
\begin{equation}
 \left (
\begin{array}{ccc}
(\Id-2\vb\wb^T ) \Jm(\xb), & -2  (\Fb(\xb), \wb)\Id, & -2\vb \Fb (\xb)^{T}\\
\mathbb{L}_1, & \Jm(\xb) - \alpha(\vb) \Id -\vb\vb^{T}(\Jm(\xb)+\Jm(\xb)^{T}), & 0 \\
\mathbb{L}_2, & -2 \wb \wb^{T} \Jm(\xb) + \wb\vb^{T}(\Jm(\xb)+\Jm(\xb)^{T}) , &\Jm(\xb)^{T}-\beta(\vb,\wb) \Id - 2\wb \vb^{T} \Jm(\xb)^{T}  \\
\end{array}
\right)\label{eqn:Jacobian2}
\end{equation}
where $\mathbb{L}_1$, $\mathbb{L}_2$ are   $n\times n$ matrices
and $\Id$ is the $n\times n$ identity matrix. To derive the above formula,  we have used the results from \eqref{eqn:alpha-beta} that 
$\nabla_{\vb}( \alpha)= \vb^{T}(\Jm^{T}  + \Jm) $,  $\nabla_{\vb}(\beta) = 2 \wb^{T}\Jm -\vb^{T}(\Jm^{T}  + \Jm) $
and 
$\nabla_{\wb}(\beta) = 2\vb^{T}\Jm^{T} $ .

In the first $n$ rows of $\tilde{\Jm}$, there are two $n\times n$
blocks which contain the term $\Fb(\xb)$ and thus vanish at the
fixed point $\xb_{s}$. So the eigenvalues of $\tilde{\Jm}(\xb_s,\vb_i,\wb_i)$ can be
obtained from the eigenvalues of its three $n\times n$ diagonal blocks: $\mathbb{N},\mathbb{M}$
and $\mathbb{K}$: 
\begin{equation*}
\begin{split}
\mathbb{N}&=(\Id-2\vb_i\wb_i^T ) \Jm(\xb_{s}), \\
\mathbb{M}&= \Jm(\xb_{s}) - \lambda_{i} \Id -\vb_{i}\vb_{i}^{T}(\Jm(\xb_{s})+\lambda_{i}\Id), \\
\mathbb{K}&= \Jm^{T}(\xb_{s})-\lambda_{i} \Id - 2\wb_{i} \vb_{i}^{T} \Jm^{T}(\xb_{s}).
\end{split}
\end{equation*}
Here  the obvious facts that $\alpha(\vb_{i})=\beta(\vb_{i},\wb_{i})=\lambda_{i}$
and $\vb_{i}^{T}\Jm^{T}=\lambda_{i}\vb_{i}^{T}$ are applied.

Now we derive the eigenvalues  of $\mathbb{N}$, $\mathbb{M}$ and $\mathbb{K}$
by constructing the corresponding eigenvectors. 
Note that $\vb_i^{T}\wb_j=\delta_{ij}$ holds under our assumption of the eigenvectors. One can verify that 
\begin{equation*}
\begin{split}
\mathbb{N} \vb_i &= 
(\Id-2\vb_i\wb_i^T)\lambda_i \vb_i = -\lambda_i \vb_i, \\
\mathbb{M}\vb_{i}&=  -2\lambda_{i} \vb_{i}\vb_{i}^{T}\vb_{i} = -2\lambda_{i}\vb_{i}, \\
\mathbb{K}\wb_{i} &= - 2\lambda_{i}\wb_{i}\vb_{i}^{T}\wb_{i}=-2\lambda_{i}\wb_{i},  \\  
\end{split}
\end{equation*}
and for all $j\neq i$,
\begin{eqnarray}
  \mathbb{N} \vb_j &=& (\Id-2\vb_i\wb_i^T )\lambda_j \vb_j = \lambda_j \vb_j,  \\
\mathbb{K}\wb_{j} &=& (\lambda_{j}- \lambda_{i})\wb_{j}- 2 \lambda_{j} \wb_{i}\vb_{i}^{T}\wb_{j}= (\lambda_{j}- \lambda_{i})\wb_{j},
\end{eqnarray}
and  with a bit more effort, 
\begin{equation}
\begin{split}
\mathbb{M}( \vb_{j}- (\vb_{i}^{T}\vb_{j}) \vb_{i})
&= \mathbb{M}\vb_{j}  - \vb_{i}^{T}\vb_{j}( \mathbb{M}\vb_{i})
=\mathbb{M}\vb_{j}  + 2\lambda_{i} (\vb_{i}^{T}\vb_{j})\vb_{i}\\
&= (\lambda_{j}- \lambda_{i})\vb_{j}  -(\lambda_{j}+\lambda_{i})  \vb_{i}(\vb_{i}^{T}\vb_{j})+ 2\lambda_{i} (\vb_{i}^{T}\vb_{j})\vb_{i}\\
&=(\lambda_{j}-\lambda_{i}) (\vb_{j} - (\vb_{i}^{T}\vb_{j})\vb_{i}).
\end{split} 
\end{equation}

Hence the eigenvalues of the Jacobian $\tilde{\Jm}$ at any fixed point
$(\xb_s,\vb_i,\wb_i)$ ($i=1,\cdots,n$)  are 
 \begin{equation}
 -2\lambda_i, \ -\lambda_i,  \ \{\lambda_j: j\neq i\}, \ \{\lambda_j-\lambda_i: j\neq
 i\}.
 \label{eqn:eig-GAD}
\end{equation}
The first and  last set of eigenvalues have multiplicity 2. The linear stability condition
is that all numbers in \eqref{eqn:eig-GAD} are negative.
Thus one fixed point $(\xb_s, \vb_{i'}, \wb_{i'})$ is linearly stable if and
only if $\lambda_{i'}>0$ and all other eigenvalues $\lambda_j<0$ for $j\neq i'$, in which case
the fixed point $\xb_{s}$ is index-$1$ saddle.


Next, we discuss some examples of GAD.

Consider first the case of a gradient system with 
$V(\xb) = \xb^T A \xb/(\xb^T \xb)$, where $A$ is a symmetric matrix.  
$V$ is nothing but the Rayleigh quotient.
A simple computation shows that the GAD for this system is given by:

\begin{equation}
\left\{
\begin{split}
\dot{\xb} = & -\frac{A \xb}{\xb^T \xb}  + \frac{\xb^T A \xb}{(\xb^T \xb)^2} \xb  +
   2 \left( \frac{\vb^T A \xb}{\xb^T \xb} - \frac{\xb^T A \xb}{(\xb^T \xb)^2} (\vb^T \xb) \right) \vb ,\\
\dot{\vb} = &  - A \vb + (\vb^T A \vb) \vb. 
\end{split}
\right. 
\end{equation}

Next, we consider an infinite dimensional example. The potential energy functional
is the Ginzburg-Landau energy for scalar fields:
$I(u) = \int_\Omega \left( \frac 12 |\nabla u|^2 + \frac 14 (u^2 - 1)^2 \right) d\xb $.
The steepest decent dynamics in this case is described by the well-known
Allen-Cahn equation:
\begin{equation}
 \partial_t u = \Delta u - (u^2-1) u.
\end{equation}
A direct calculation gives the GAD in this case:
\begin{equation}
\left\{
\begin{split}
 \partial_t u  & = \Delta u - (u^2-1) u - 2(\Delta u -(u^2-1) u, v) v, \\
 \partial_t v  & = \Delta v - ( 3u^2-1) v - (\Delta v -(3u^2-1) v, v) v ,
\end{split}\right.
\end{equation}
where the inner product is defined to be:
\begin{equation*}
 (u, v) = \int_\Omega u(\xb) v(\xb) d\xb.
\end{equation*}

Clearly both the SDD and the GAD depend on the choice of the metric, the inner product.
If we use instead the $H^{-1}$ metric,  then the SDD becomes the Cahn-Hilliard equation
and the GAD changes accordingly.

\section{High index saddle points}

GAD can also be extended to the case of finding high index saddle points. 
We will discuss how to
generalize it to index-$2$ saddle points here.  There are
two possibilities: Either the Jacobian $\Jm$ at the saddle point has
one pair of conjugate complex eigenvalues or it has  two real eigenvalues at
the saddle point.
We discuss each separately.

%
%
Intuitively, the picture is as follows.
We need to find the projection of the flow, $\Fb(\xb)$,
on the tangent plane, say $P$, of the two dimensional unstable manifold of the saddle point,
and change the direction of the flow on that tangent plane.
For this purpose, we need to find the vectors $\vb_{1}$ and $\vb_{2}$ that span $P$. 
In the first case, we assume that  the unstable  eigenvalues at the saddle point
are $\lambda_{1,2}=\lambda_{R}\pm i \lambda_{I}$. 
In this case there are no real eigenvectors corresponding to $\lambda_{1,2}$.
However, for any vector $\vb$ in $P$, $(\nabla \Fb)\vb$ simply rotates  $\vb$ inside  $P$.
Hence, $\vb_{2}$ can be taken as 
$(\nabla \Fb) \vb_{1}$ if we have already found some $\vb_{1}\in P$.
The latter can be accomplished using the original dynamics in \eqref{GAD-non-grad}.

To see how one should modify the flow $\Fb$ on the tangent plane, we write
$$\Fb=c_{1}\vb_{1}+c_{2}\vb_{2} + \sum_{j>2}c_{j}\vb_{j}.$$
Using the fact that the eigen-plane of $(\nabla \Fb)^{T}$ corresponding to 
$\lambda_{R}\pm i \lambda_{I}$, which is
spanned by $\wb_{1}$ and $\wb_{2}=(\nabla \Fb)^{T} \wb_{1}$, 
is orthogonal  to $\vb_{j}$ for all $j>2$, we can derive a linear system for $c_1$ and
$c_2$ by taking the 
inner product of ${\Fb}$ and $\wb_{1}$,$\wb_{2}$.
The solution of that linear system is given by:
\begin{equation}
c_{1} = \frac{a_{22}f_{1}-a_{12}f_{2}}{a_{11}a_{22}-a_{21}a_{21}}, \ c_{2} = 
\frac{a_{11}f_{2}-a_{21}f_{1}}{a_{11}a_{22}-a_{21}a_{21}}
\label{eqn:c2}
\end{equation}
where $a_{ij}=(\wb_{i},\vb_{j})$ and $f_{j}=(\Fb(\xb),\wb_{j})$ for $i,j=1,2$.
 The \gad for the $\xb$ component is  $$\tilde{\Fb}=\Fb -2c_{1}\vb_{1}-2c_{2}\vb_{2}.$$
{To s}ummarize, we obtain the following dynamical system:
\begin{equation}
\left\{
\begin{split}
\dot{\xb} &= \Fb -2c_{1}\vb_{1}-2c_{2}\vb_{2}, \\
\dot{\vb}_{1} &= (\nabla \Fb(\xb))\vb_{1} - \alpha(\vb_{1})\vb_{1},\\
\dot{\wb}_{1} &= (\nabla \Fb(\xb))^T \wb_{1} -\beta(\vb_{1},\wb_{1})\wb_{1},\\
\vb_{2} &= \nabla \Fb(\xb)\vb_{1} ,\\
\wb_{2} &= (\nabla \Fb(\xb))^{T}\wb_{1},
\end{split}
\right.
\label{eqn:GAD-ng-c2}
\end{equation}
where $c_1,c_2$ are given by  \eqref{eqn:c2} and 
$\alpha,\beta$ are defined by \eqref{eqn:alpha-beta}.

  If the Jacobian has two positive real eigenvalues
at the saddle point, say, $\lambda_{1}>\lambda_{2}>0 \geq \lambda_{3} >  \cdots$,
let us define a new matrix by the method of deflation:
\begin{equation}
\Jm_{2} := \nabla \Fb- \frac{ (\vb_{1}, (\nabla \Fb) \vb_{1})}{(\vb_{1},\vb_{1}) (\wb_{1},\vb_{1})} \vb_{1}\wb_{1}^{T}.
\label{eqn:J2-def}
\end{equation}
It is not difficult to see that if $\vb_{1}$ is an eigenvector of $\nabla \Fb$ 
corresponding to $\lambda_{1}$, then ${\Jm}_{2}$ shares the same eigenvectors as $\Jm$, and 
 the eigenvalues of $\Jm_{2}$ become $0,\lambda_{2},\lambda_{3},\cdots$.
The largest eigenvalue of $\Jm_{2}$ at the index-$2$ saddle point becomes $\lambda_{2}$. 
One can then use the dynamics \eqref{GAD-non-grad-v} associated with the new matrix $\Jm_{2}$ to find $\vb_{2}$. Therefore, 
we obtain the following  index-$2$ GAD
\begin{equation}
\left\{
\begin{split}
\dot{\xb} &= \Fb -2c_{1}\vb_{1}-2c_{2}\vb_{2}, \\
\dot{\vb}_{1} &= (\nabla \Fb(\xb))\vb_{1} - \alpha_{1}\vb_{1},\\
\dot{\wb}_{1} &= (\nabla \Fb(\xb))^T \wb_{1} - \beta_{1}\wb_{1},\\
\dot{\vb}_{2} &= \Jm_{2} \vb_{2}  - \alpha_{2}\vb_{2} ,\\
\dot{\wb}_{2} &= \Jm_{2}^{T} \vb_{2}  - \beta_{2}\wb_{2},
\end{split} 
\right.
\label{eqn:GAD-ng-r2}
\end{equation}
with the initial normalization condition $(\vb_{1},\vb_{1}) =(\vb_{2},\vb_{2}) =(\wb_{1},\vb_{1}) =(\wb_{2},\vb_{2})=1$.
$c_1$ and $c_2$ are given in the same way as shown above   \eqref{eqn:c2} and 
$\alpha_{1,2},\beta_{1,2}$ are defined as follows to enforce that the 
normalization condition is preserved : 
$\alpha_{1} =(\vb_{1}, (\nabla \Fb(\xb))\vb_{1}), 
\beta_{1}= 2(\wb_{1}, (\nabla \Fb(\xb))\vb_{1})-\alpha_{1}$ 
and 
$\alpha_{2} =(\vb_{2},\Jm_{2}\vb_{2}), 
\beta_{2}= 2(\wb_{2}, \Jm_{2}\vb_{2})-\alpha_{2}$.

The generalization to higher index saddle points with real eigenvalues is obvious.

\section{Examples}

\subsection{Analysis of a gradient system}

To better understand the dynamics of GAD, let us consider the case when 
a different relaxation  parameter is used
for the direction $\vb$:
$$\left\{
\begin{split}
 \dot{\xb}  &= - \nabla V(\xb) + 2(\nabla V, \vb) \vb,  \\
 \tau \dot{\vb}  &= - \nabla^2 V(\xb)\vb + (\vb, \nabla^2 V \vb) \vb.  
\end{split}
\right.$$
To simplify the discussions, we consider the limit as $\tau \to 0 $. 
In this case, we obtain a closed system for $\xb$:
\begin{equation}
 \dot{\xb}  = - \nabla V(\xb) + 2(\nabla V, \vb(\xb)) \vb(\xb), 
 \label{eqn:gad-x}
 \end{equation}
where $\vb(\xb)$ is the eigenvector of $\nabla^{2}V(\xb)$ associated with the smallest eigenvalue. 
Now we consider the following two dimensional system:
$$V(x,y)=\frac14(x^{2}-1)^{2}+\frac12 \mu y^{2}$$
where $\mu$ is a positive parameter.  $\xb_{\pm}=(\pm 1,0)$ are two stable fixed points
and $(0,0)$ is the index-$1$ saddle point.
The eigenvalues and eigenvectors of the Hessian at a point $\xb=(x,y)$ are
\begin{align*}
\lambda_{1}&=3x^{2}-1 \text{ and } \vb_{1}=(1,0),\\
\lambda_{2}&=\mu \text{ and } \vb_{2}=(0,1).
\end{align*}
Therefore, the eigendirection picked by GAD  is 
\begin{equation}
\begin{cases}
\vb_{GAD}(\xb) = \vb_{1}, &\mbox{ if }  |x| < \sqrt{\frac{1+\mu}{3}},\\
\vb_{GAD}(\xb) = \vb_{2}, &\mbox{ if }  |x| > \sqrt{\frac{1+\mu}{3}}.   
\end{cases}
\end{equation}
Consequently, by defining 
$$V_{1}(x,y)=-\frac14(x^{2}-1)^{2}+\frac12 \mu y^{2},$$
and
$$V_{2}(x,y)=\frac14(x^{2}-1)^{2}-\frac12 \mu y^{2},$$
we can write the \gad \eqref{eqn:gad-x} in the form of a gradient system driven by the new potential:
\begin{equation}
\displaystyle
V_{GAD}(\xb)= V_{1}(\xb)\cdot 1_{|x|<\sqrt{\frac{1+\mu}{3}}} (\xb)
+V_{2}(\xb)\cdot 1_{|x|>\sqrt{\frac{1+\mu}{3}}} (\xb)
\label{eqn:Vgad}
\end{equation}
where $1_{\cdot}(\xb)$ is the indicator function.
Note that $V_{GA\com{D}}$ is \textit{not} continuous at the lines
$x=\pm \sqrt{\frac{1+\mu}{3}}$. The point $(0,0)$ becomes the unique local minimum
of $V_{1}$, with the basin of attraction  $\{(x,y):-1<x<1\}$.  
Outside of this basin of attraction, the flow goes to
$(x=\pm\infty,y=0)$ and the potential $V_{1}$ falls to $-\infty$ .
For $V_{2}$, the point $(0,0)$ is the unique local maximum and 
all solutions go to $(x=\pm 1, y=\pm \infty)$.

\begin{figure}[htbp]
\begin{center}
\includegraphics[width=\textwidth]{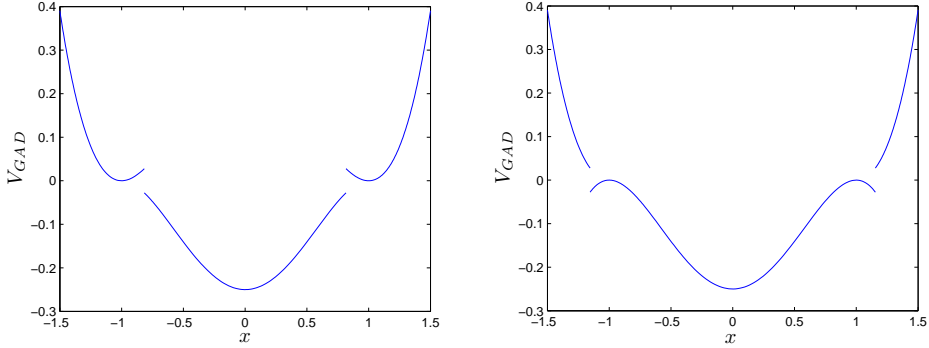}
\caption{The discontinuity of $V_{GAD}(x,y=0)$ at the location $x=\pm \sqrt{\frac{1+\mu}{3}}$.
Left: $\mu<2$; Right: $\mu>2$.}
\label{fig:Vgad_dis}
\end{center}
\end{figure}

\begin{figure}[htbp!]
\begin{center}
\includegraphics[width=0.86\textwidth, height=0.32\textwidth]{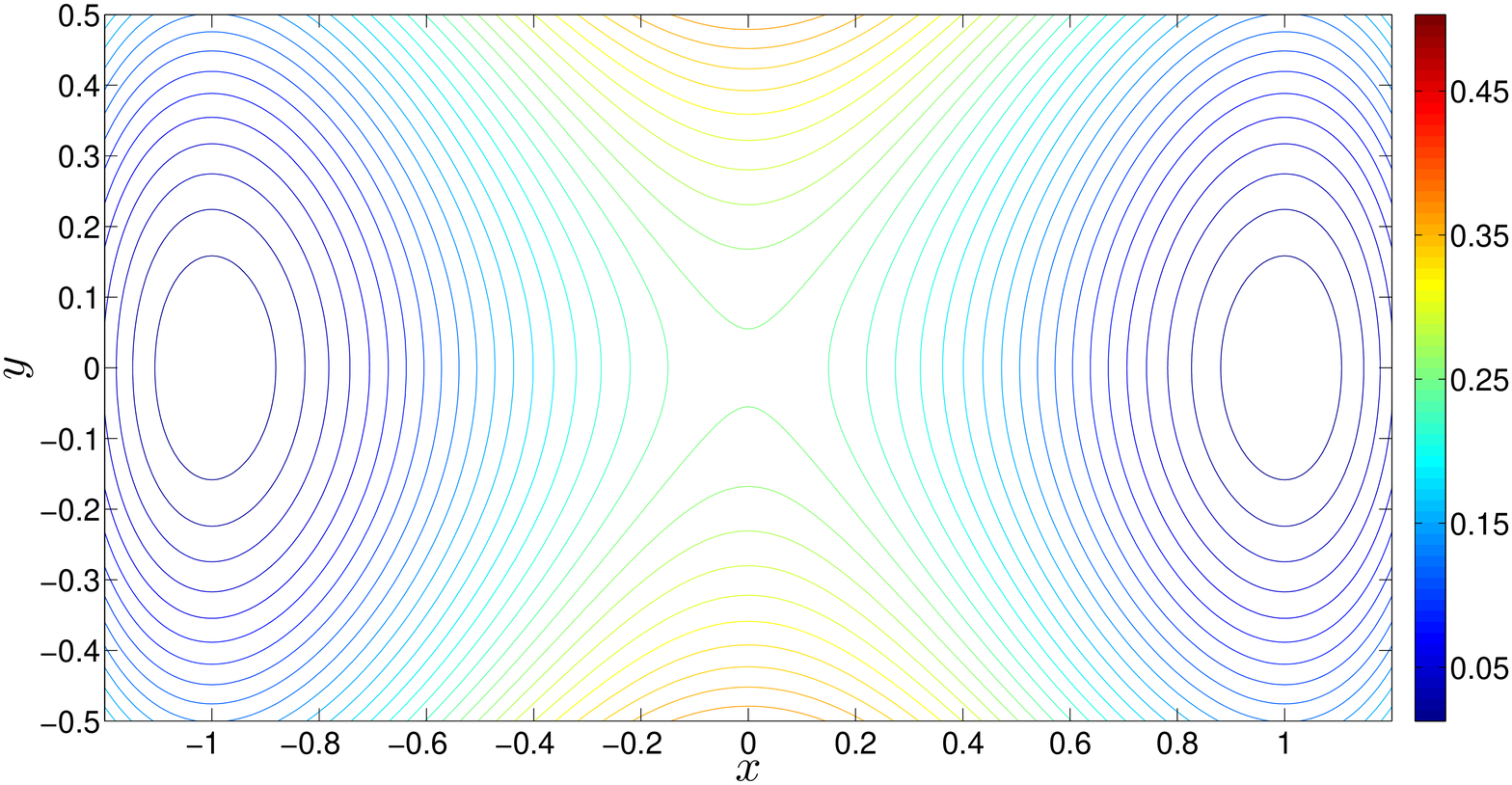}
\includegraphics[width=0.86\textwidth, height=0.32\textwidth]{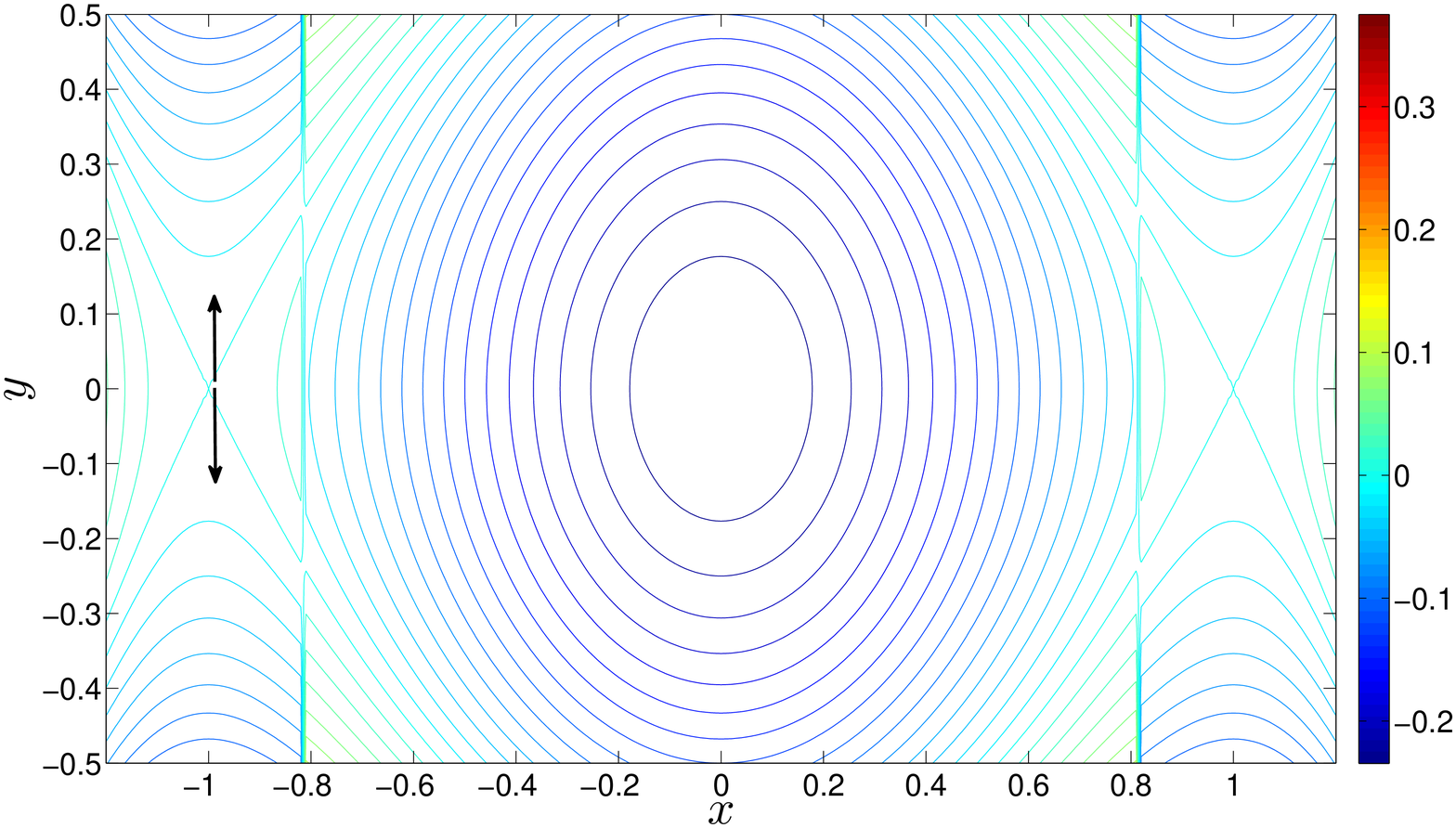}
\includegraphics[width=0.86\textwidth, height=0.32\textwidth]{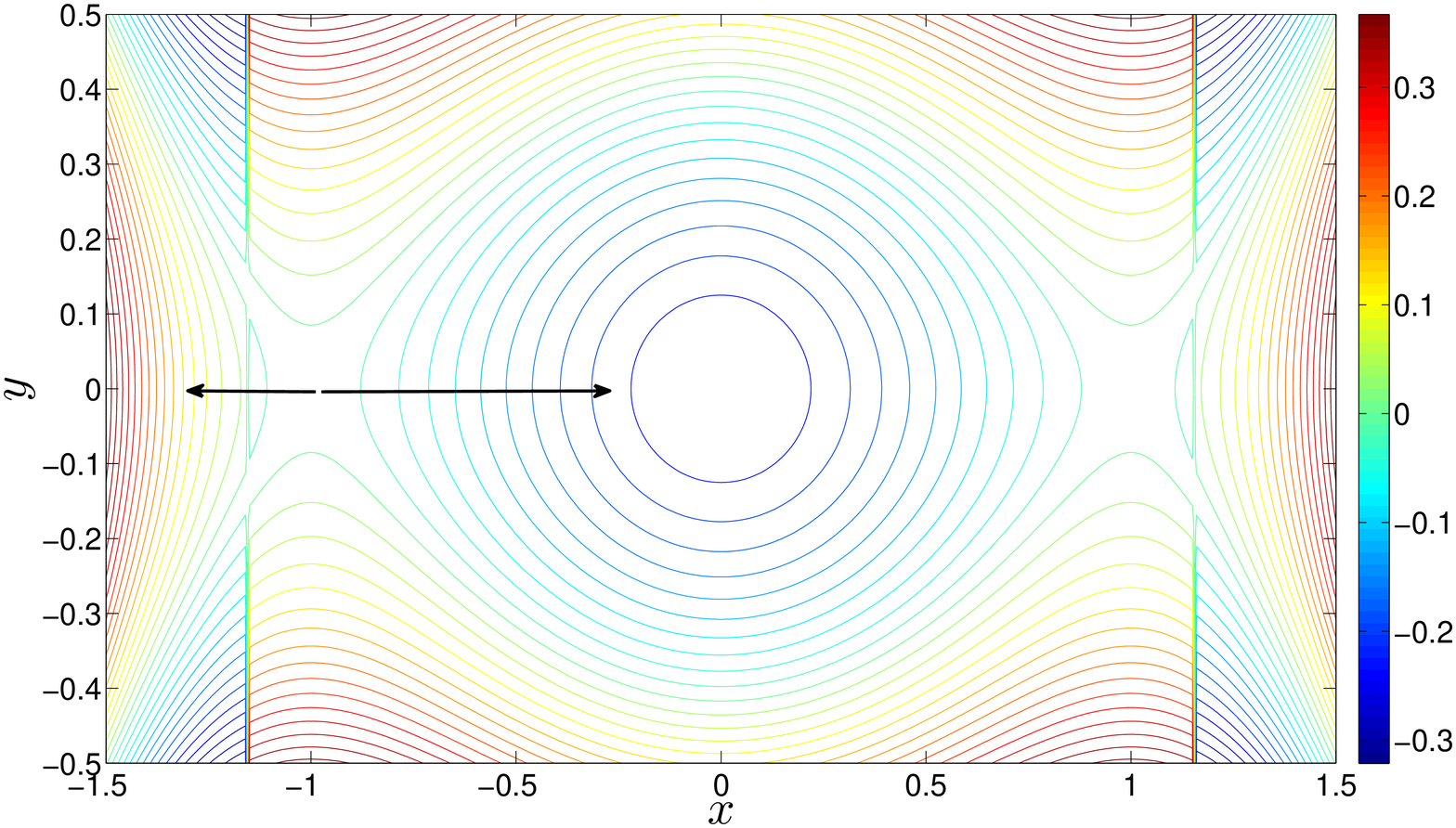}
\caption{The contour plots of $V$, $V_{GAD}$ for $\mu=1$ and $V_{GAD}$ for $\mu=3$, from the top to the bottom, respectively. For the plot of $V_{GAD}$,  $V_{1}$ lies in the middle  region
$-\sqrt{\frac{1+\mu}{3}}<x<\sqrt{\frac{1+\mu}{3}}$  and the 
 $V_{2}$ lies at the two sides.
The arrows show the flow directions of the \gad  \eqref{eqn:gad-x}.}
\label{fig:grad}
\end{center}
\end{figure}

If we start the \gad with the initial value $\xb_{\pm} = (\pm 1,0)$, 
then there are two different situations according to whether $\mu>2$ or $\mu<2$. 
Although $\xb_{\pm}$ becomes a saddle point for any $\mu\neq 2$,
the unstable direction for $\mu<2$ is $\pm \vb_{2}$
while the unstable direction for $\mu>2$ is $\pm \vb_{1}$,
as illustrated in figure \ref{fig:grad}. Furthermore, from figure \ref{fig:grad}
and the above discussion, 
it is clear that the basin of attraction of the point $(0,0)$ associated with the potential $V_{GAD}$
is the region $-\sqrt{\frac{1+\mu}{3}}<x<\sqrt{\frac{1+\mu}{3}}$ 
for $\mu<2$ and   $-1<x<1$ for $\mu>2$. (which 
is larger than the basin of attraction for the Newton-Raphson method, confirmed by numerical calculation.)  Consequently, the GAD with an initial value 
$(x_{0},y_{0})$  \emph{near}  the 
local minimum $\xb_{\pm}$ of $V$ converges to the point $(0,0)$ of our interest
when $\mu>2$ and $|x_{0}|<1$.

This discuss suggests that GAD may not necessarily converge globally and instabilities
can occur when GAD is used as a numerical algorithm.
When instabilities do occur, one may simply reinitialize the initial position 
or the direction.

\subsection{Lorenz system} 

\begin{figure}[htbp!]
\begin{center}
\includegraphics[width=0.8\textwidth]{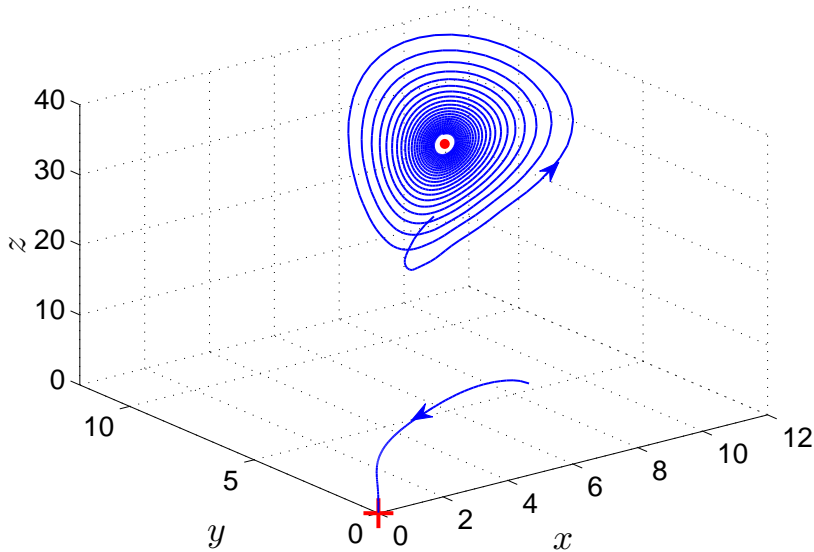}
\caption{The trajectories of GAD for the Lorenz system starting from two initial points. They converge
to the index-2 saddle point $Q_{+}$ (marked by the dot) and the index-1 saddle point $O$ (marked by ``$+$'') 
respectively.}
\label{fig:lorenz-1}
\end{center}
\end{figure}

\begin{figure}[htbp!]
\begin{center}
\includegraphics[width=0.8\textwidth]{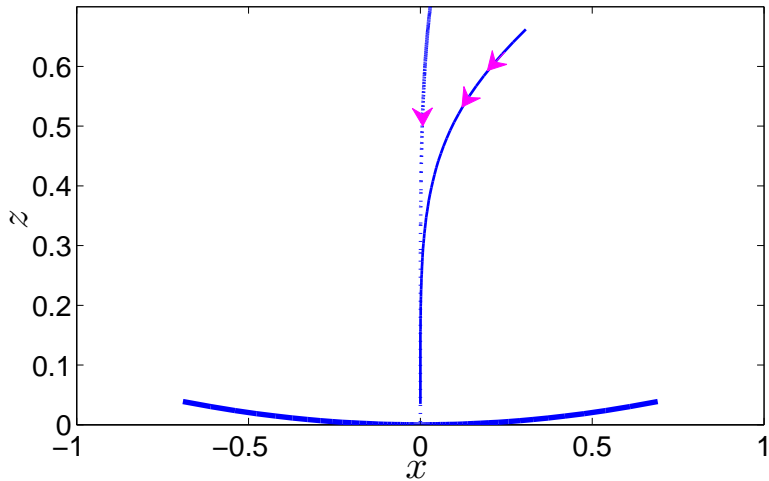}
\caption{How the GAD trajectories approaches the saddle point $O$. The curve with two arrows
is the trajectory of index-$1$ GAD for the Lorenz system; the curve with single arrow is the trajectory of 
index-$2$ GAD for the time reversed Lorenz system. 
The unstable manifold of $O$, which is tangent to the $z=0$ plane, is also shown.}
\label{fig:lorenz-2}
\end{center}
\end{figure}

Consider
\begin{equation} \label{eqn:Lorenz}\left\{\begin{array}{lcl}
\dot{x} &=& \sigma (y-x), \\ \dot{y} &=& \rho x - y -xz,\\
\dot{z} &=& -\beta z + xy. \end{array}\right.
\end{equation}
The parameters we use are 
$\sigma=10$, $\beta=\frac83$ and $\beta=30$.   There are three fixed points:
the origin $O=(0,0,0)$  and 
two symmetric  fixed points
$$Q_{\pm}=(\pm \sqrt{\beta(\rho-1)},\pm
\sqrt{\beta(\rho-1)},\rho-1).$$
$O$ is an index-$1$ saddle point. 
The Jacobian at $Q_{\pm}$ has one pair of complex conjugate eigenvalues
with positive real part.
In our calculation, we prepare the
initial directions $\vb_{0}$ and $\wb_{0}$ by running the GAD 
for long time starting from random initial conditions for
$\vb$ and $\wb$ while keeping $\xb$ fixed,
although this is not entirely necessary.
Figure ~\ref{fig:lorenz-1} shows two solutions of GAD.
For the index-$1$ saddle point $O$, figure \ref{fig:lorenz-2} depicts how 
the trajectory of GAD converges to it. 
It can be seen that the component of the original force $\Fb$ along the unstable 
direction of $O$ is nearly projected out, thus the trajectory will not be affected by 
the unstable flow in that direction and avoids departing the saddle point. 
Therefore the trajectory tends to follow the stable manifold toward the saddle point when the trajectory is close enough
to the saddle point. Similar behavior is seen for the case of searching
 the point $Q_{+}$ which has one pair of complex eigenvalues.
The trajectory surrounding  $Q_{+}$ in the figure \ref{fig:lorenz-1} 
spirals to $Q_{+}$ and these spirals are 
closer and closer to the unstable manifold of $Q_{+}$ in the original Lorenz dynamics, 
which looks like a  twisted disk. 
The convergence rate of the spiraling trajectories in GAD is very slow
because the real part of the complex eigenvalues 
($\lambda=0.1474 \pm 10.5243\text{ i} $) in the original dynamics is rather small 
compared with its imaginary part.

If we reverse time $t\to -t$, we have the time-reversed Lorenz system, 
in which the origin $O$  becomes an index-$2$ saddle point. 
We can apply the index-2 GAD algorithm \eqref{eqn:GAD-ng-r2} to search for this saddle point.
The GAD trajectory in this case is also plotted in the figure \ref{fig:lorenz-2}.  
It is similar to the situation of GAD applied to the original Lorenz system
in the sense that the GAD trajectory nearly follows the $z$ axis when approaching
 the limit point $O$.
Indeed, as far as the $\xb$-component is concerned, the linearized \gad for the original
Lorenz system and the time-reversed one are the same.
From the proof of the Proposition (particularly, note that the eigenvalues of $\mathbb{N}$
are $-\lambda_{i}$ and $\lambda_{j}$), 
it is not hard to see that the eigenvalues of the linearized \gad at the point $O$ are all negative and 
 have the same absolute values as the eigenvalues of the 
original dynamics, and the two 
dynamics  share the same eigenvectors (again, we mean the  $\xb$ component of the GAD). 
Thus, since the 
change $t\to -t$ does not change the absolute values of the eigenvalues of the original dynamics, the \gad for the original
and time reversed Lorenz system have the same eigenvalues: 
$\lambda_{1} = -23.3955$, $\lambda_{2}=-2.6667$, $\lambda_{3}=-12.3955$. 
The two linearized GAD flows near the point $O$ are the same: $\xb(t)=e^{-23.3955t}\vb_{1}+e^{-2.6667t}\vb_{2}+e^{-12.3995t}\vb_{3}$,
where $\vb_{1,2,3}$ are the eigenvectors:  $\vb_{2}=(0,0,1)$, and $\vb_{1}$, $\vb_{3}$
are in the $z=0$ plane.   As $t\to +\infty$, 
we then have $\xb(t) \sim e^{-2.6667t}\vb_{2}$. 
This explains why both trajectories in the figure \ref{fig:lorenz-2}  
follow the $z$ axis when approaching the saddle point $O$.

\subsection{A PDE example with nucleation} 

\begin{figure}[htbp!]
\begin{center}
 \includegraphics[width=0.85\textwidth]{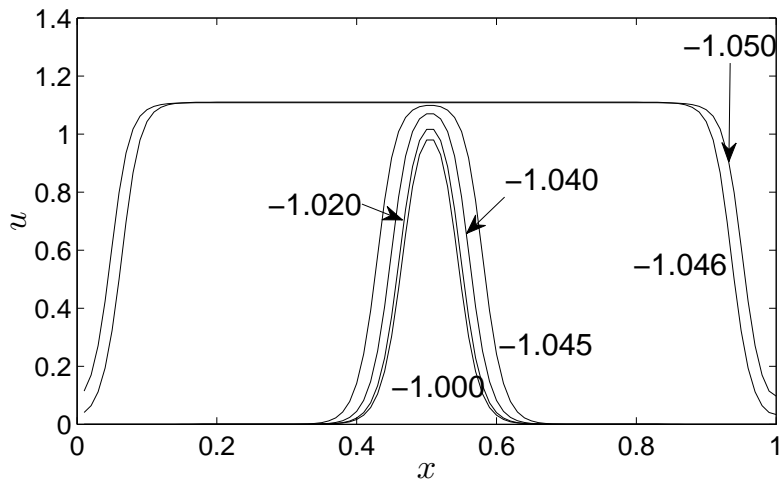}
\caption{The profiles of saddle points of the   example
\eqref{eqn:ngre1} ($\delta=0.01$). Only the component $u$ is shown since $v=\frac12 u^{2}$ at the saddle point.
 From inside to outside, the
values of $\mu$ are $-1.000$, $-1.020$, $-1.040$, $-1.045$, $-1.046$,
$-1.050$.} \label{fig:case0_saddleprofiles}
\end{center}
\end{figure}

Let us consider the following reaction-diffusion system 
on the domain $x\in[0,1]$ with periodic boundary condition:
\begin{equation}
\left\{
\begin{split}
\frac{\partial u}{\partial t} &= \delta \Delta u +  {\delta}^{-1}f(u,v),\\
 \frac{\partial v}{\partial t} &= \delta \Delta v + {\delta}^{-1}g(u,v),
\end{split}\label{eqn:RD2}
\right.
\end{equation}
where
\begin{equation}
\label{eqn:ngre1}
\begin{cases}
 f(u,v) &= (u-u^3+1.2)v+ \frac12 \mu u,\\
 g(u,v) &=  \frac12 u^2 -v.
\end{cases}
\end{equation}
The parameter $\delta$ is fixed at $0.01$ and we allow the parameter  $\mu$ to vary.
There are two stable (spatially homogeneous) solutions for certain range of $\mu$: $\ub_{+}=(u_{+},v_{+})$ and $\mathbf{0}=(0,0)$.
If one uses the square-pulse shape function as a initial guess in the Newton-Raphson method,
no convergence can be achieved in most situations. 
We applied the index-$1$ GAD method to this example.
The initial conditions for GAD are constructed by adding a small amount of 
perturbations around either stable solutions: $\ub_{+}$ or $\mathbf{0}$. 
We observed that for a fixed value of $\mu$,  the solutions of GAD constructed this way 
converge to the same saddle point.
The different saddle points obtained from GAD at different values of $\mu$ are plotted
in figure \ref{fig:case0_saddleprofiles}.
It is also numerically confirmed that these saddle points indeed have index $1$ and the unstable manifold
goes to $\ub_{+}$ in one unstable direction and to $\mathbf{0}$ in the opposite
unstable direction.
It is interesting to observe the dependence of the saddle point on the parameter $\mu$
and that such a dependence is highly sensitive when $\mu$ is close to $-1.046\sim-1.045$.
In fact, there exists a critical value $\mu^{*}$ in this narrow interval at which the spatially extended 
system \eqref{eqn:RD2} has a subcritical bifurcation, which does not appear in the corresponding ODE system
without spatial dependence. 
We refer to \cite{subcrit2010} for further discussions about this point.

\section{Concluding remarks}

We expect that
GAD is particularly useful for handling  high dimensional system in the sense that  
it should have a larger basin of attraction for finding  saddle points, than, for example,
the Newton-Raphson method.
There are many questions one can ask about GAD.  
One question is the convergence of GAD as time goes to infinity.
Our preliminary result shows that GAD does not have to converge.
For finite dimensional systems, there is always local convergence
near the saddle point. The situation for infinite dimensional systems, i.e. PDEs,
seems to be much more subtle.
Another interesting point is whether one can accelerate GAD.
For the problem of finding local minima, 
many numerical algorithms have been proposed and they promise to have much faster 
convergence than SDD.  It is natural to ask whether analogous ideas
can also be found for saddle points.

{\bf Acknowledgement:}  The work presented here was supported in part
by AFOSR grant FA9550-08-1-0433.
The authors are grateful to Weiguo Gao and Haijun Yu and the second referee for helpful discussions.

\bibliography{./ngre_copy}

\bibliographystyle{siam} 

\end{document}